
\magnification=1200
\overfullrule=0mm
\nopagenumbers
\vskip 1cm
\centerline {\bf Polynomial Interpolation and}
\centerline{\bf a Multivariate Analog of 
Fundamental Theorem of Algebra}
\bigskip
\centerline {\bf H. Hakopian and M. Tonoyan}
\bigskip

Starting with univariate polynomial interpolation we arrive to a natural generalization of fundamental theorem of algebra in form of a system of multivariate algebraic equations. We  characterize  the case when the algebraic system has maximal number of distinct  solutions. This gives also a characterization of  correct systems of knots for multivariate Lagrange interpolation. We present a conjecture for the case of solutions with multiplicities, i.e. for the case of multivariate  Hermite interpolation. 

Next we discuss the above generalization from the point  of view of  differential equations. According to the generalization we  develop a multivariate  setting which includes natural  analogs of well-known results  concerning the first-order normal  system of ordinary differential equations and the nth order normal differential equations.

\bye